# A model for the optimal design of a supply chain network driven by stochastic fluctuations

K. Petridis[1], A. K. Chattopadhyay[2] and P. K. Dey[1]

[1]Aston Business School, Aston University, Aston Triangle,
Birmingham B4 7ET, UK

&

[2]Nonlinearity and Complexity Research Group, Aston University,
Aston Triangle, Birmingham B4 7ET, UK

**ABSTRACT**

Supply chain optimization schemes have more often than not underplayed the role of inherent stochastic fluctuations in the associated variables. The present article focuses on the associated reengagement and correlated renormalization of supply chain predictions now with the inclusion of stochasticity induced fluctuations in the structure. With a processing production plant in mind that involves stochastically varying production and transportation costs both from the site to the plant as well as from the plant to the customer base, this article proves that the producer may benefit through better outlay in the form of higher sale prices with lowered optimized production costs only through a suitable selective choice of producers whose production cost probability density function abides a Pareto distribution. Lower the Pareto exponent, better is the supply chain prediction for cost optimization. On the other hand, other symmetric (normal) and asymmetric (lognormal) distributions lead to upscaled costs both in terms of inlays and outlays. While this is an averaged out statistics over large time regimes, transient features may still affect such probabilistic predictions and offset results. The predictions are shown to be in good harmony with model results shown.

*Keywords*: Optimal supply chain network design, MINLP, Facility Location Problems, Stochastic Models, Noise

# 1. Introduction

An inviolable aspect of business organization is the distribution of supply lines both at the input side of the business as well on its output deliverables, together with the supply chain management (SCM) of its throughput overall. Due to the rapid economic globalization, majority of the operations conducted from manufacturing to transportation avenues, and from warehousing to the customer base, are conducted by supply chain contractors or third party logistics (3PL) companies. Recent research in logistics developments shows that in the years to come, approximately 80% of economic transactions will be based on services. Thus, the better the design of the supply chain operations, the better the service level the customers will experience as nowadays in the majority of products transported and sold throughout the world, customers are not brand loyal, thus any stockout may result in sales reduction and future loss of income for firms. (Gruen, Corsten, and Bharadwaj 2002; Gruen and Corsten 2004).

While the subject has hugely benefitted from paradigmatic studies in the realm of supply chain theories involving deterministic variation of associated variables and parameters, very little has been done in connection with the impact of stochastic perturbations in probabilistically predicting the qualitative and quantitative assessment from the supply chain model. Optimizing the design of a supply chain provides an "ideal" image of the real situation that, by construction, is not amenable to conventional mathematical modeling (Mohammadi Bidhandi et al. 2009; Seferlis and Giannelos 2004). The issue here is the randomized nature of the data produced from the supply chain performance profile that are mathematically categorized as stochastic in design. There are of course stochastic programming models, where each unique scenario is associated with a corresponding probability of occurrence, yet such models have historically failed to incorporate explicit stochastic effects in their formulation. The real situation entails mismatches in the operations conducted among the nodes of the supply chain, in such a way that affects the levels of upstream and downstream decision (Tamas 2000).

Identifying stochasticity in the operations is not enough, though, to provide real and stochastic decisions (Santoso et al. 2005). The mismatches that take place in the operations of supply chain network design are taken into consideration with fluctuations in the variables and the question that is raised is how different types of fluctuations affect the supply chain network design.

The primary objective of this paper is to develop a mathematical model of a supply chain that accounts for all inherent stochastic fluctuations of the system and its parameters.

In this paper, a multi-stage multi-echelon model is presented that hierarchically incorporates functional interactions between plants, warehouses, customer zones, and thereby multiple echelons in turn. The first and final links of the supply chain are considered to be fixed and only the quantities of products produced (for plants) and transported (for customers) are provided. A graphical representation of the supply chain network is provided in Figure 1.

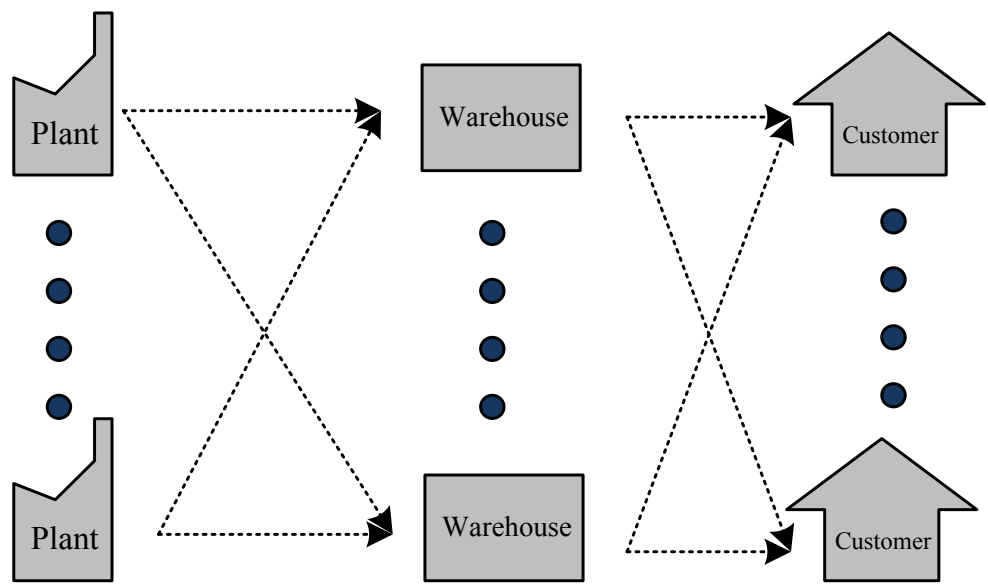

**Figure 1** Multi-stage, multi-echelon supply chain network

The warehouses are assumed not to be installed. Thus, based on the intermediate link, the supply chain network is constructed. The proposed model has two stages; in the first stage the optimal design of the supply chain network is calculated based on stochastic demand assuming that is normally distributed Petridis (2013); the second stage is fed with levels of decisions from the first stage in order to compute any shortfalls in demand, and to compute the expected lead time as well as the quantities that must be produced in order to vanish the stock out instances.

The stock out instances, which are defined as the absolute difference between demand and the quantities of products delivered to customer, are divided into two categories based on a threshold decided by the decision maker (DM). In many cases, a stock out instance may not just affect the service level and therefore the perception of the customers towards a specific product but may lead to penalty costs due to a contract clause (especially in the food supply chain industry). Due to the fact that besides holding inventory, warehouses are assumed also to serve as production facilities, and that the inventory can serve as raw materials in order to cover the deficits in demand. The magnitude of production quantities is assigned to a corresponding production cost that is added to the total cost function of the 1st stage. Due to stock out instances, the expected lead time (*ELD*) keeps increasing. The fact that warehouses are used as production facilities in the supply chain reduce the expected lead time but may increase the overall cost significantly, leading to a trade-off between cost and service quality. The aforementioned procedure is graphically represented in Figure 2.

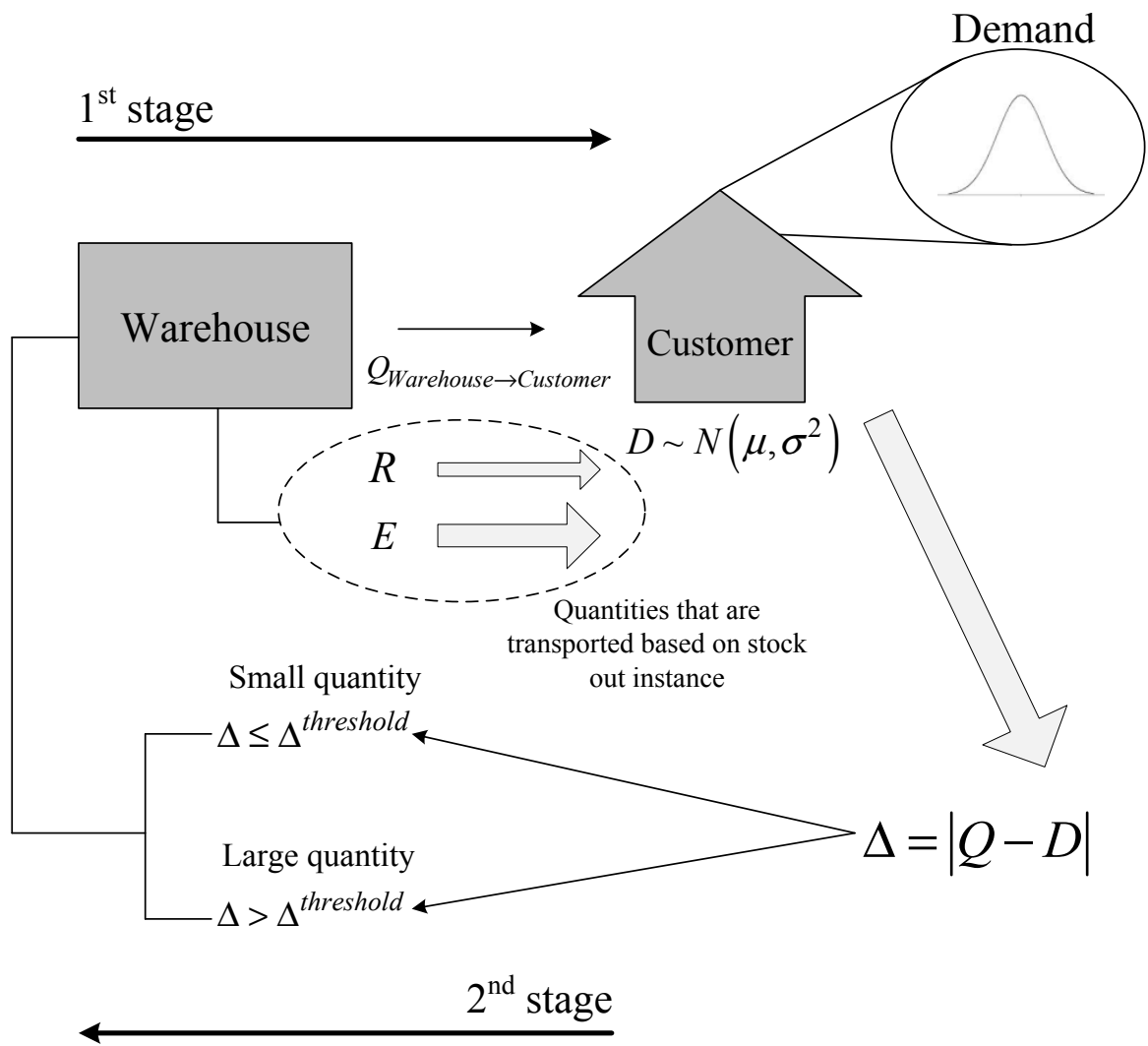

**Figure 2** The two-stage supply chain model

In the proposed two stage model, there is a new characteristic that of the incorporation of stochastic noise into the design of the supply chain network. Assuming different types of noise representations in terms of their respective probability density functions (Gaussian, Lognormal, and Pareto), the supply chain model is analyzed to quantify which of these PDFs ensure cost minimization through optimization rationale perpetrated across the entire supply chain network. Fluctuations are not formulated through different parameter distribution representations, but directly are introduced to the examined variables. Through the proposed approach the following research questions are derived:

- The fluctuation of which distribution approach better the "ideal" situation?
- The fluctuation of which distribution increases cost?

## 2. Literature survey

The optimal supply chain network design (SCND) problem has been extensively examined in the literature. The majority of the models used are taken from mathematical programming disciplines and are roughly divided into two categories; a) to steady state models and b) to multi-period models. Due to the absence of time, the decisions of continuous variables represent average levels per specified unit of time Beamon (1998).

In their approach Tsao and Lu (2012) proposed a Nonlinear Programming (NLP) model providing an integrating framework for the facility location and inventory allocation problem with cost discounts. A two-phase approximation approach was deployed as a solution to provide numerical results that could demonstrate the impact of different simulated data to the supply chain decisions and cost. In a novel development directed towards fluctuation incorporation, Tsiakis, Shah, and Pantelides (2001) used a Mixed Integer Linear Programming (MILP) model where both binary and continuous variables were considered with the objective of assigning uncertainty in the structure of the hierarchical variables e.g. demand as deterministic uncertainties in their respective numbers, without explicit incorporation of statistical stochastic terms; the first are used for network representation while the latter for facility capacity and flows of goods throughout the channels of the supply chain network Melo, Nickel, and Saldanha da Gama (2006). Similar studies have been proposed considering the demand uncertainty measuring the customers' service level through the calculation of lead time and normally distributed demand You and Grossmann (2008)

Petridis (2013). The formulation of an agile or flexible supply chain network with the use of a heuristic algorithm as a solution procedure was also proposed by (Pan and Nagi 2010) as a means of incorporating certain non-deterministic fluctuations perpetrating changed functionalities in the supply chain.

Close-loop supply chains (CLSC) are generally used to model the reusability and recycling of products (ICT, food etc). In a recent work, (Jindal and Sangwan 2014) proposed a fuzzy MILP model in order to capture the uncertainty in demand, cost and other parameters. Similar studies have been proposed in the literature using mathematical programming techniques for the optimal close loop supply chain network design (CLSCND) (Harold Krikke, Bloemhof-Ruwaard, and Van Wassenhove 2001); (H. Krikke, Bloemhof-Ruwaard, and Van Wassenhove 2003).

Recent works focus mostly on biomass based supply chain networks due to a global turn toward bioenergy production. In their work Grigoroudis, Petridis, and Arabatzis (2014) proposed a Data Envelopment Analysis (DEA) based algorithm for the optimal biomass supply chain network design. The optimal design of a forest supply chain network was proposed by (Arabatzis et al. 2013). In this work, they employed a Lagrangian Relaxation algorithm Fisher (1985) to design the fuel-wood supply chain, considering demand uncertainty. The optimal design of biofuel supply chain network has been also examined using a Monte Carlo simulation approach to provide a sensitivity analysis for various parameters (Kim, Realff, and Lee 2011).

The use of multiple objective functions provides more realistic approaches to real world problems. In such domains, multi-objective programming (MOP) models have been traditionally employed, including the optimal design of chemical supply chains (Guillén-Gosálbez and Grossmann 2009), biofuel/biomass supply chains (Gebreslassie, Yao, and You 2012), (Liu and Papageorgiou 2013), in forest supply chains (Arabatzis, Petridis, and Kougioulis 2014) or considering green supply chains with environmental factors (Wang, Lai, and Shi 2011).

The introduction of noise realization has been examined in many production-allocation systems (including the supply chain network design problem). The main modeling method for noise representation is optimal control. (Riddalls and Bennett 2002) have proposed a multi-echelon control model in order to describe a production-allocation supply chain network. In their work, they assumed noise-corrupted demand and system delays. A game theoretic, Stackelberg game model was proposed by (An et al. 2007), where through a collaborative approach, a noise (read fluctuation) reduction scheme was propounded. Noise, in terms of uncertainty, has also been modeled through different demand and supply scenarios identifying disruptions production process (Baghalian, Rezapour, and Farahani 2013). A decision support system was proposed by (Acar, Kadipasaoglu, and Schipperijn 2009) where the performance of service level or customer satisfaction was examined through a simulation based study. Uncertainty has been modeled by adding noise to demand parameter or by sampling from statistical distributions the data regarding other operations like supply. Although more realistic, explicit incorporation of multiplicative noise based routines have rarely come across in this literature, partly due to computational difficulty as also due to the minimalist nature of most problems considered.

The majority of the works proposed in the supply chain literature are modeling uncertainty or stochasticity by performing a sensitivity analysis in the parameters of the models (bounds, demand, supply etc.). However, in many scenarios, the parameters of a supply chain or adding an index to each variable do not represent the misalignments of the conducted operations in the supply chain. The information mismatch is also another phenomenon that is not addressed adequately within the frameworks of the aforementioned types of stochastic models. In this work, the problem of designing in an optimal way the supply chain network design is performed incorporating different types of noise in the variables of the study.

The present work lays the foundation of a new optimization routine in which the noise representations are taken from well known statistical distributions and then the cost optimization, based on a cost function, is done by optimizing both with respect to the stochastic variables as also with respect to the stochastic parameters.

## 3. Methods
### 3.2 Model concepts
Here the model concepts of the study are presented. There are two types of concepts: deterministic and stochastic. In the first case, the model is solved in the deterministic limit without any noise addition in the variables while in the second, different uncertainty representations are modeled with varying noise distributions (normal, lognormal and Pareto).

### 3.2.1 Deterministic

The optimization model presented here provides levels of decisions for the quantities of a single product, even though extensions can be also considered. The model presented here is an extension of MINLP model of Petridis (2013).

For better understanding of the optimization model, each node (stage) is assigned with an index. The first stage (plants) is denoted by *i*, the second (warehouses) by *j*, and the third by *k*. In the following context, the constraints of the problem are presented. As each plant has a limited capability given its resources, raw material etc., the production capability of each plant *i* is upper and lower bounded.

$$P_i \leq P_i^U \quad \forall i \tag{1}$$

$$P_i \geq P_i^L \quad \forall i \tag{2}$$

From (1) and (2), $P_i^U$ and $P_i^L$ model production upper and lower bounds which are assumed to be known a priori.

As the produced quantities travel downstream (from the production to the end customer), the following mass balance constraint is considered, modeling the fact that the produced quantities by plant *i* should equal to the quantities transported from plant *i* to warehouse *j*.

$$P_i = \sum_j Q_{ij} \quad \forall i \tag{3}$$

The quantities entering warehouse node should be equal to those that exit that node (from warehouse *j* to customer zone *k*).

$$\sum_i Q_{ij} + I_j = \sum_k Q_{jk} \quad \forall j \tag{4}$$

Finally the quantities transported throughout the supply chain end to the customers' end should be greater than or equal to the demand of each customer. The demand is assumed to follow a statistical distribution that is already known.

$$\sum_j Q_{jk} \geq D_k \quad \forall k \tag{5}$$

As mentioned above, the warehouse facilities are not known a priori and are decided after the optimization model. Generally, decisions of yes or no type are introduced with binary variables. The connection between the plant $i$ and warehouse $j$ and customer $k$ and warehouse $j$, is assumed to exist only if warehouse $j$ exists.

$$X_{ij} \leq Y_j \quad \forall i, j \tag{6}$$

$$X_{jk} \leq Y_j \quad \forall j, k \tag{7}$$

In constraints (6) and (7) variables binary $Y_j$ model whether warehouse $j$ will be installed in position $j$ or not and $X_{1,2}$ to model the connection between nodes 1 and 2.

The quantities are transported from one node to another only if the corresponding connection exists.

$$Q_{ij} \leq Q_{ij}^U \cdot X_{ij} \quad \forall i, j \tag{8}$$

$$Q_{jk} \leq Q_{jk}^U \cdot X_{jk} \quad \forall j, k \tag{9}$$

Finally, the warehousing quantities can be computed through the following constraints (a Lagrange multiplier type approach):

$$W_j \geq a_j \cdot \sum_i \left(Q_{ij} + I_j\right) \quad \forall j \tag{10}$$

$$W_j \leq W_j^U \cdot Y_j \quad \forall j \tag{11}$$

The objective function of the model is to minimize the overall cost:

$$\begin{aligned} TC = \sum_i c_i^P \cdot P_i + \sum_i \sum_j c_{ij}^{VTR} \cdot Q_{ij} + \sum_i \sum_j c_{ij}^{FTR} \cdot X_{ij} \\ \sum_j \sum_k c_{jk}^{FTR} \cdot Q_{jk} + \sum_j \sum_k c_{jk}^{VTR} \cdot X_{jk} + \sum_j c_j^{IN} \cdot Y_j \end{aligned} \tag{12}$$

In objective function (12), the first term represents the production cost, the second and fourth terms respectively represent variable transportation costs, while the third and fifth terms account for the fixed variable costs and the final term represents the installation (or capital) cost.

As there may be shortfalls in demand (unsatisfied demand) due to bad information or scheduling, natural disasters that may disrupt this chain etc., the following parameter is introduced:

$$\Delta_k = \left| D_k - \sum_j Q^*_{jk} \right| \ \forall k \tag{13}$$

This parameter $\Delta_k$ serves as a critical threshold to define the level of unforeseen expenses to be expected in emergent conditions. This threshold is divided into the following ranges: $\left[\Delta^L, \Delta^M\right]$ and $\left(\Delta^M, \Delta^U\right]$; the threshold is divided into a low values range of supply insufficiency as indicated by the first range, and in a high values of supply insufficiency, corresponding to very large shortfalls in providing the demanded quantity. In the examined ranges is stands that $\Delta^L = \min_k \{\Delta_\kappa\}$, $\Delta^U = \max_k \{\Delta_\kappa\}$ while $\Delta^M = 0.5 \cdot \left(\Delta^L + \Delta^U\right)$.

In equation (13), the absolute difference between the demand representation of each customer *k* and the transported quantities computed from the previous stage (stage 1) is shown. In order to model the magnitude of the failure in customer's satisfaction, the following constraints are now introduced:

$$\Delta^L \cdot \lambda_k \leq Q^U_k \leq \Delta^M \cdot \lambda_k \ \forall k \tag{14}$$

$$\Delta^M \cdot \zeta_k \leq Q^O_k \leq \Delta^U \cdot \zeta_k \ \forall k \tag{15}$$

$$\lambda_k + \zeta_k = 1 \ \forall k \tag{16}$$

Constraints (14) and (15) are introduced to model the deficit in demand under or over a pre specified threshold. Binary variables $\lambda_k$ and $\zeta_k$ are mutually exclusive as any shortfall in demand can be characterized as over or under a specific threshold but cannot fall in both categories, as indicated in (16).

In stock out instances, several corrective actions should be undertaken to improve the service level without significantly increasing the cost. In case where the deficit in demand belongs to the interval above the predetermined threshold (new quantities are produced from the fixed inventory kept in warehouse *j* ($I_j$) along with the corresponding deficit of demand for this particular customer *k*.

If the quantity that corresponds to low values of supply insufficiency, binary variable that corresponds to the low range becomes 1 for some indices of customers of *k* ($k_1$), then $\lambda_{k_1} = 1$ and based on the predetermined range in (14), warehouses which are also assumed to be production plants holding inventory used for manufacturing purposes, will have to produce additional quantity equal to $R_{jk}$ as seen in (18). Constraint (18) provides a value that corresponds to the quantity to be produced based on constraint (17), as binary variable $K_{jk}$ takes a value of 1 if-f $\lambda_k$ equals to 1. In that case, the quantity that will be eventually produced by warehouse $j$ in order to facilitate a medium stockout occurred at customer $k$ should be more than $\gamma_{jk} \cdot H_{jk}$; $\gamma_{jk}$ stands for the production coefficient of warehouse warehouse $j$ for each customer $k$ and $H_{jk}$ is a minimum level of inventory stored for the production of the necessary quantity in warehouse $j$ in order to facilitate a medium stockout occurred at customer $k$. Constraint (19) models the occurrence of a large stockout instance while the production quantity that is needed to be sent to customer $k$ from warehouse $j$ is defined as $E_{jk}$ and should be more than a warehouse's $j$ production rate ($\beta_{jk}$) multiplied by the sum of the overall inventory held at warehouse $j$ and stockout occurred in customer $k$ as in (20).

$$K_{jk} \leq \lambda_k \quad \forall k \tag{17}$$

$$R_{jk} \geq \gamma_{jk} \cdot H_{jk} \cdot K_{jk} \quad \forall j,k \tag{18}$$

$$\Omega_{jk} \leq \zeta_k \quad \forall k \tag{19}$$

$$E_{jk} \geq \beta_{jk} \cdot (I_j + \Delta_k) \cdot \Omega_{jk} \quad \forall j,k \tag{20}$$

Constraint (5) is now reformulated as follows:

$$\sum_j Q_{jk} + Q_k^U + Q_k^O = D_k \quad \forall k \tag{21}$$

Based on (21), the probability of stock out or over stock can be computed. If demand is normally distributed $D_k \sim N(\mu, \sigma^2)$, then $\Delta_k$ is assumed to be normally distributed as well (as difference of two random variables that are normally distributed). Thus, stock out and overstocking probabilities are introduced with the following constraints:

$$P_k(Q^U) = \frac{1}{2}\left[1 + erf\left(\frac{Q^U - \bar{\Delta}_k}{\sigma\sqrt{2}}\right)\right] \quad \forall k \tag{22}$$

$$P_k(Q^O) = 1 - P_k(Q^U) \quad \forall k \tag{23}$$

The service level in supply chain can be easily quantified with the expected lead time, namely, the amount of time needed for a product to be delivered to the customer after order placement.

The Expected Lead Time (ELD) is computed based on the following equality:

$$ELD_k = T^u \cdot P_k\left(Q^O\right) \cdot \zeta_k + T^l \cdot P_k\left(Q^U\right) \cdot \lambda_k \quad \forall k \tag{24}$$

This leads to the new objective function as follows:

$$TC1 = TC + \sum_j \sum_k c_{jk}^{PO} E_{jk} + \sum_j \sum_k c_{jk}^{PU} \cdot R_{jk} + \sum_k \sigma \cdot \sqrt{ELD_k} \tag{25}$$

In the objective function presented in (25), σ is the standard deviation of unsatisfied demand for customer $k$ such that: $\sigma = \dfrac{\sqrt{\sum_k \left(\Delta_k - \bar{\Delta}_k\right)^2}}{n-1}$ and $\bar{\Delta}_k$ the mean unsatisfied demand.

From the above analysis, the following levels of decision are derived from each stage:

- 1$^{st}$ stage
    - Produced and transported quantities
    - Selected warehouses and capacity
    - Supply chain network
    - Demand deficit
- 2$^{nd}$ stage
    - Stockout and overstocking probabilities
    - Expected Lead Time (*ELD*)
    - Quantities that should be produced to cover unsatisfied demand

### 3.2.2. Stochastic

In stage 2, the new variables introduced in the model with the addition of noise are introduced with the following equalities:

$$\tilde{P}_i = P_i + \eta \quad \forall i \tag{26}$$

$$\tilde{Q}_{ij} = Q_{ij} + \eta \quad \forall i, j \tag{27}$$

$$\tilde{Q}_{jk} = Q_{jk} + \eta \quad \forall j, k \tag{28}$$

```
for e = 1....n
   for e′ = 1....n
      X̃_{ee′} = X_{ee′} + η_{ee′}                                              (29)
   end
end
```

In equations (26-28), the basic variables of the study are replaced by the corresponding variables with noise driven representation ($\eta$). The procedure based upon each noise representation is incorporated in the variables is shown in (29).

### 3.2 Implementation of Deterministic and Stochastic model

In this section the implementation of the deterministic and stochastic models are presented in Figures 3 and 4. In the first case, the model is solved deterministically. As seen in Figure 3, initially the MIP problem is solved while the shortfalls in the demand are computed for each customer. In order to measure the magnitude of stock out instances, it is assumed that if $\Delta_k$ is more than the average stock out quantities, there is a large deficit in meeting demand and thus the expected lead time for demand satisfaction will be larger than in the case where this deficit is not of that magnitude. The final step as seen in Figure 3, is the calculation of the MINLP model where the expected lead time, the probabilities of over and under stocking instances, and the levels of variables are provided.

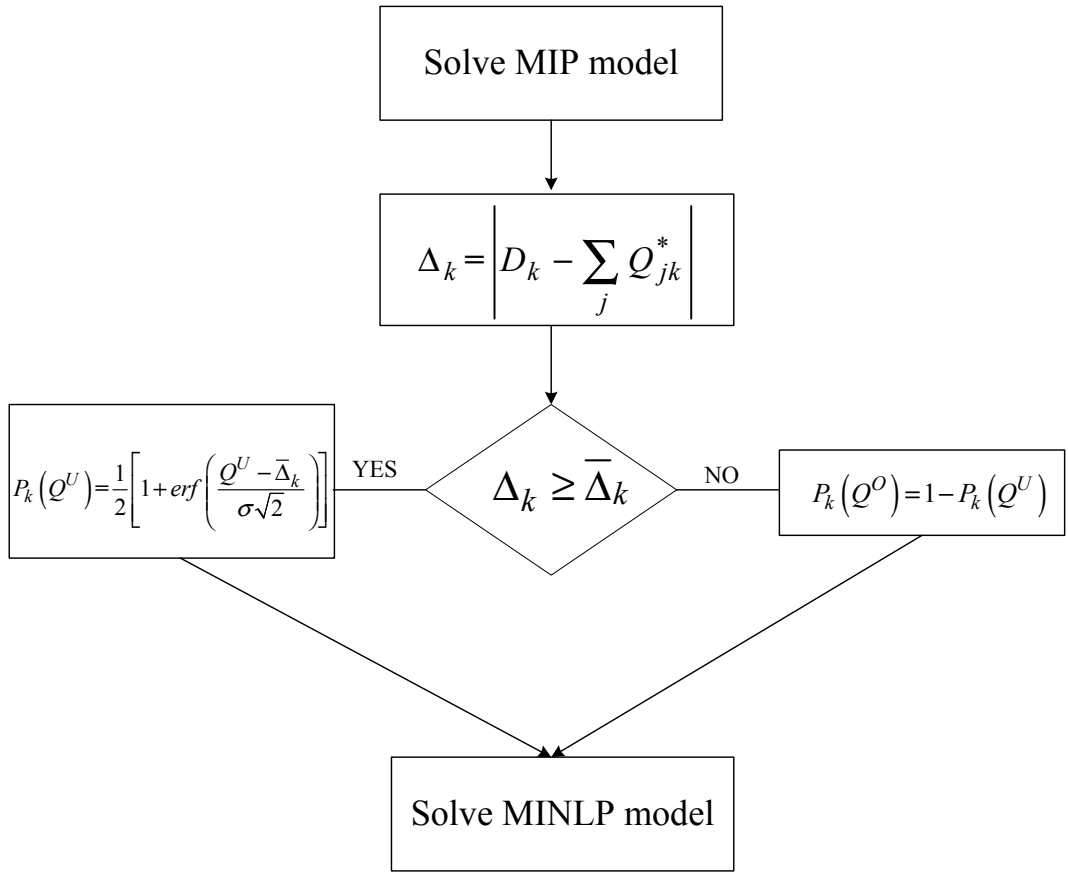

**Figure 3** Flowchart of the deterministic implementation

On the contrary, in the stochastic case as presented in the second stage of the analysis, the MINLP model is solved for different noise representations for the basic variables that concern the production and transportation of flows as described in equations (26-28). The introduction of noise into the variables is implemented using procedure (29). For each new variable, the MINLP may yield a feasible solution (optimal, local optimal of integer) or an infeasible solution. A counter is introduced to model each instance where MINLP model yields a feasible solution.

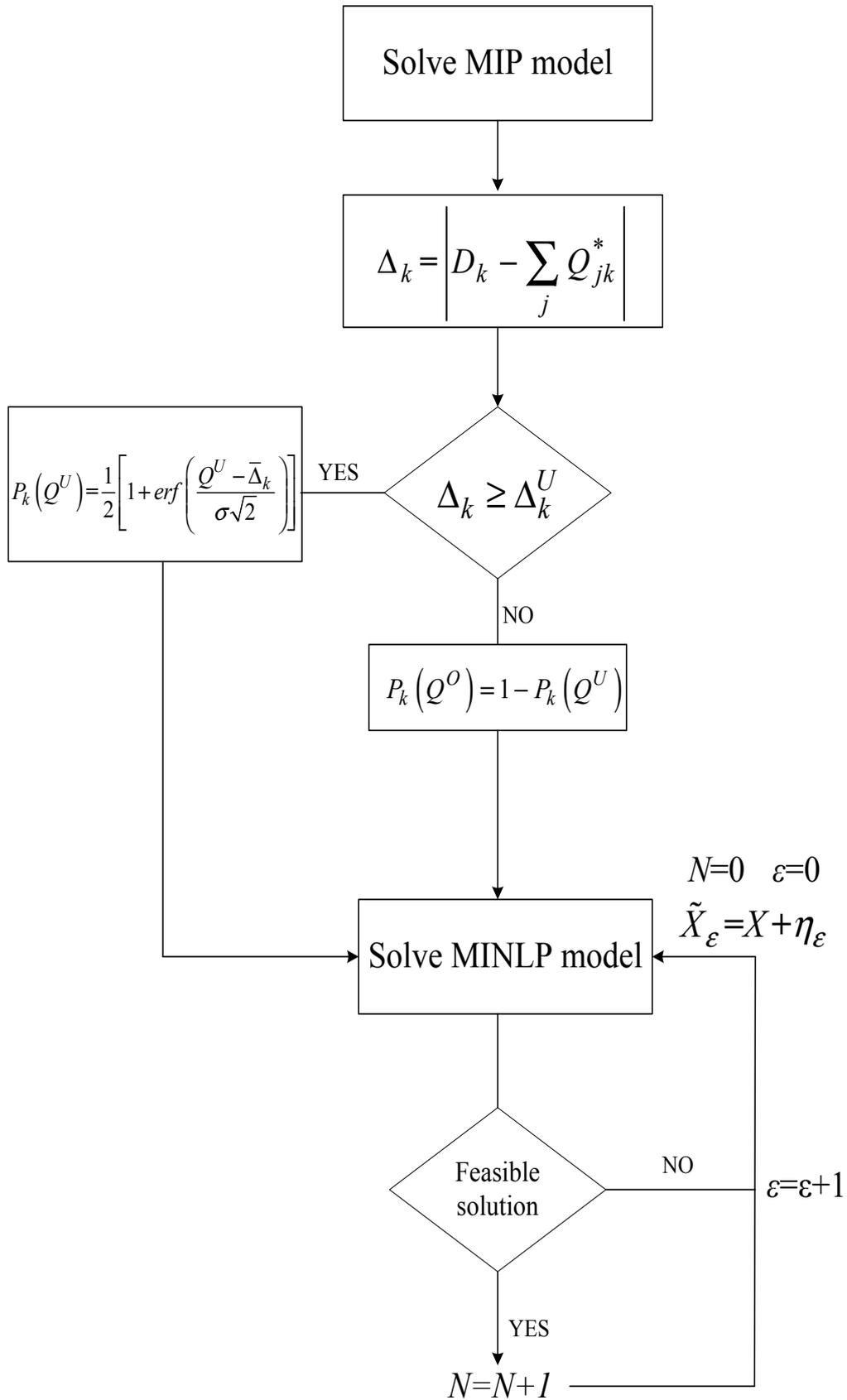

**Figure 4** Flowchart of the stochastic implementation

After the end of the loops for the noise representations, the average of noise representations are computed based on the following averaging scheme:

$$\overline{X}_e = \frac{1}{N} \sum_{e'} X_{ee'} \qquad (30)$$

Based on equation (30), the Root Mean Square (RMS) variable is defined through the following formula:

$$X^{RMS} = \sqrt{\frac{1}{C_2^N} \sum_{\substack{e_1, e_2 \in E \\ e_1 \neq e_2}} \overline{X}_{e_1} \cdot \overline{X}_{e_2}} \qquad (31)$$

In the calculation of RMS from equation (31), $C_2^N$ stands for the set of two combinations of set $N$ where $N \geq 2$ and $N$ is the feasible solution instances.

### 3.3 Risk assessment

In this section the different noise representations are presented. The representations of different distributions of noise used are presented below and graphically illustrated in Figure 5:

- Gaussian noise
- Lognormal noise
- Pareto noise for various alpha levels namely:
    - $\alpha=0.01$
    - $\alpha=0.05$
    - $\alpha=0.5$
    - $\alpha=0.99$

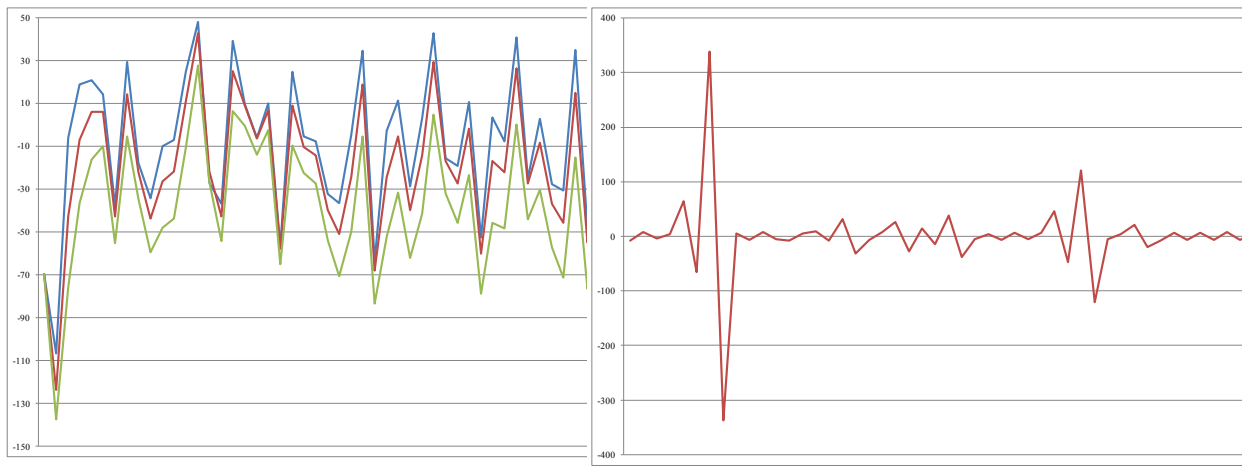

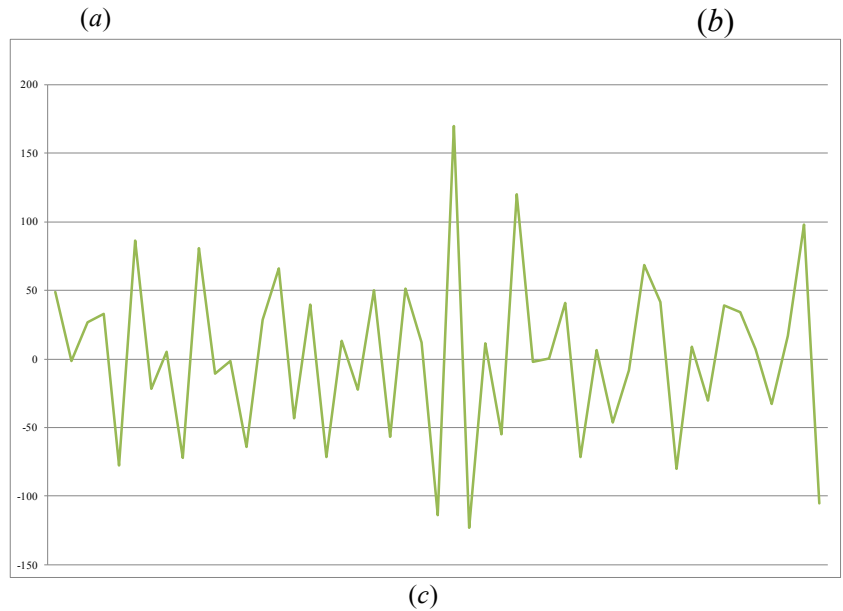

**Figure 5** Noise representations: a) Pareto noise for α=0.01 (blue line), α=0.5 (red line), α=0.99 (green line), b) Lognormal noise, c) Gaussian noise.

## 4. Results

The results of the study are focused on the results of the variables with noise (stochastic model) compared to those without noise (deterministic model). Due to its dimension, the variable that models the production of each plant *i* is presented for all noise representations. The variables $Q_{ij}$ and $Q_{jk}$ concern the transported quantities from plant *i* to warehouse *j* and from warehouse *j* to customer *k*. Due to the reasonably large number of nodes considered at each echelon, in abeyance of the statistical measure, heatmap plots are presented. In our estimation, we have assumed this nodal number to be 20 $(|I|=|J|=|K|=20)$, In Figure 6, the results of variable $P_i$ for different noise realizations are shown. Our stochastic model shows that the optimized cost with a Pareto noise distribution for *Pareto exponent* a→0 comes closest to the deterministic prediction while Pareto distributions with larger exponent values as well as other distributions, like lognormal or Gaussian, lead to poor cost optimization schemes.

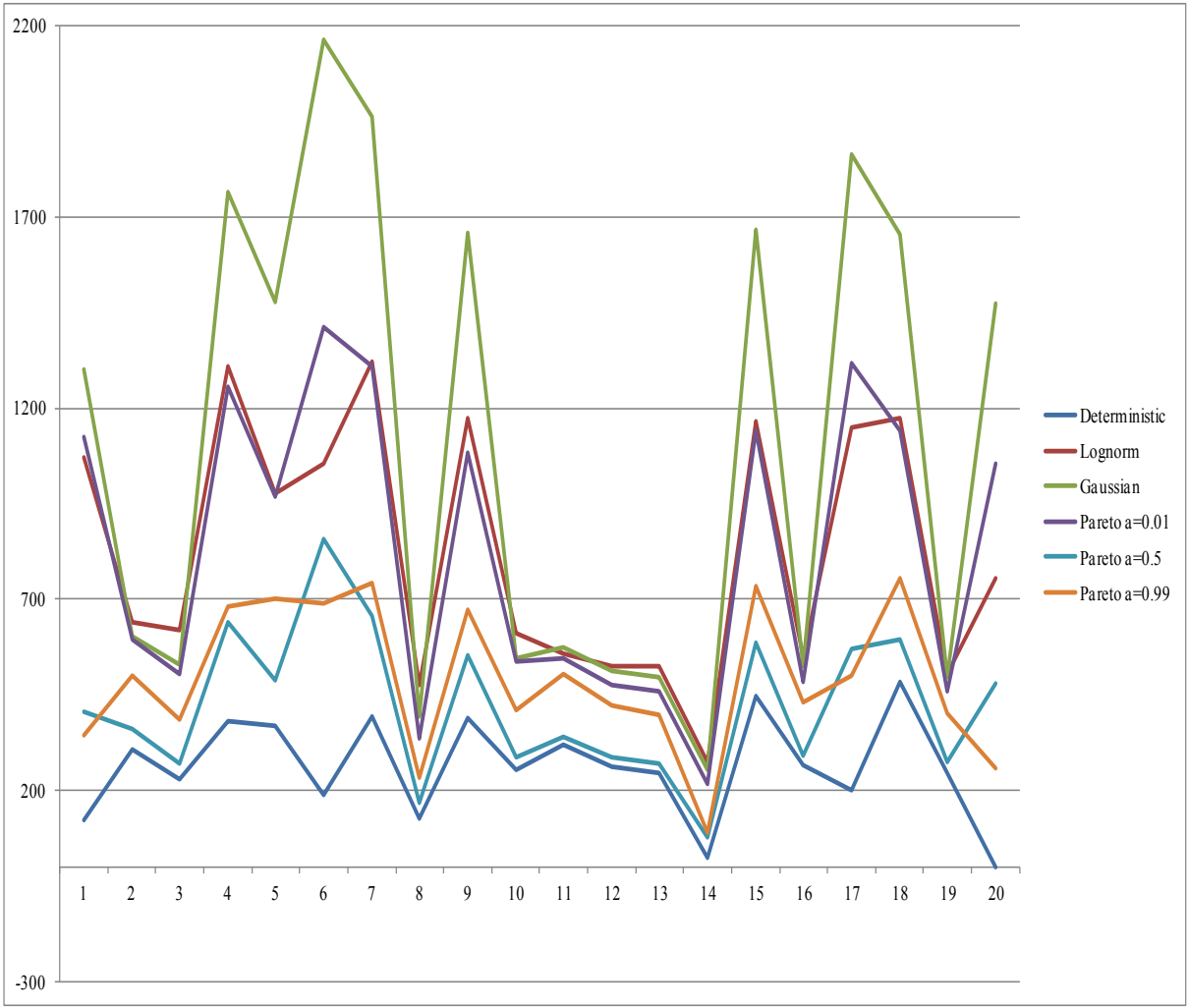

**Figure 6** Results for variable $P_i$ for different noise representations compared to the deterministic

In Figure 7, the results of the difference of deterministic variable $Q_{ij}^{det}$ with $Q_{ij}^{noise}$ which correspond to the results of transported quantities from plant *i* to warehouse *j*, are presented. If the results are closer to 0, then it can be concluded that the addition of that specific noise does not has a significant impact to the overall Supply Chain Network design. In Figure 7a), it can be seen that Pareto noise with *a*=0.01 that most of the area lies in the range of $[-100, 100]$ (green and purple color). This means that the fluctuations from the deterministic values of variable $Q_{ij}$ may be from -100 to 100. From Figure 7b), it can be seen that the fluctuations are increasing to the range $[-400, 200]$, yet in this case the "bumps" increase. In Figure 7c), it can be seen that although most the area lies in the range of $[-500, 500]$ the "bumps" reduce based on the previous case. In Figure 7d) it can be seen that there are a lot of fluctuations most of which lie in the range of $[-100, 100]$. Finally in Figure 7e) there is approximately the same image as in Figure 7d), nevertheless the majority of the area lies on the range of $[-100, 0]$.

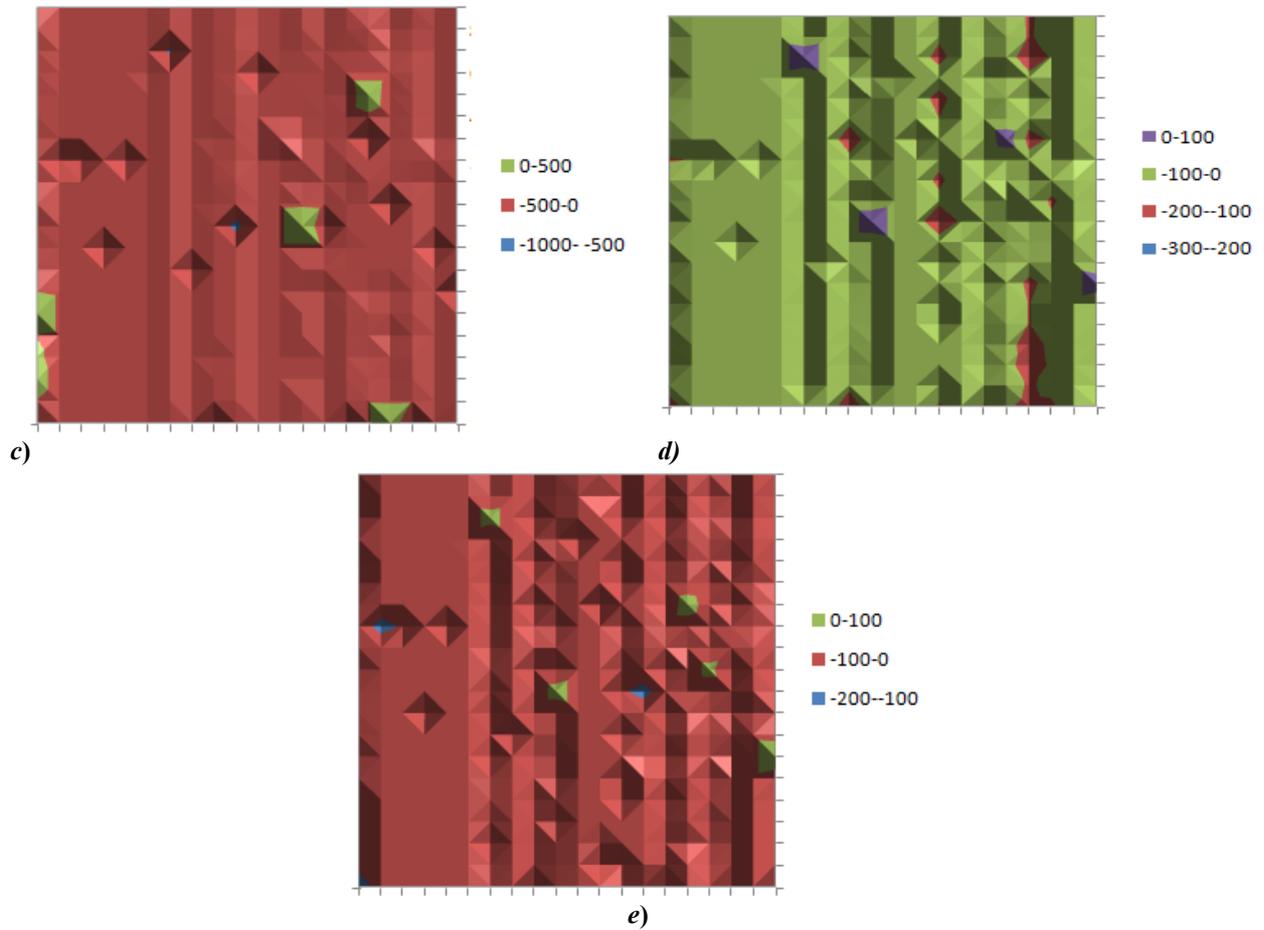

*c)* *d)*

*e)*

**Figure 7:** Heatmaps for the difference of deterministic values of $Q_{ij} - Q^*_{ij}$ for a) Pareto noise with $a$=0.01, b) Pareto noise with $a$=0.5, c) Pareto noise for $a$=0.99, d) Gaussian noise, e) Lognormal noise.

In order to provide an overall measure of the results presented above, the standard deviation of the difference of the results derived after the introduction of each noise representation with the deterministic ones is presented in the following table (Table 1).

**Table 1:** Standard deviation of the difference of deterministic value of variables with noise representation

| Noise Representation $\left(Q_{ij}^{det} - Q_{ij}^{noise}\right)$ | Standard Deviation $(\sigma)$ |
|---|---|
| Pareto $a$=0.01 | 27.92 |
| Pareto $a$=0.5 | 97.65 |
| Pareto $a$=0.99 | 70.16 |
| Gaussian Noise | 29.39 |
| Lognormal Noise | 18.97 |

## 5. Discussion and conclusions

The optimal design of a supply chain network may be oriented from the customer's perspective, namely "pull" systems or from production's orientation, namely "push" systems (Spearman and Zazanis 1992). In "pull" systems, demand drives production remains conserved while in the second case, production is fixed based on demand estimation. However, in most cases the optimal design of supply chain network is constructed around parameter values that approximate the upper and lower bounds of transported quantities.

In several studies, stochasticity has been introduced either as different scenarios or integrating a statistical distribution into the parameter (expected value), in order to capture the characteristic of uncertainty (Santoso et al. 2005). None of these approaches, though, reflect the real situation as the uncertainty is measured on the parameter and not on the variable, aside of the fact that such subroutines can only lead to implicit uncertainty measures at best and inaccurate predictions at worst. The previous statement can be easily understood with the following example. Assuming that the well-known 'bullwhip' effect occurs (Lee, Padmanabhan, and Whang 2004); the variable that corresponds to the quantities that are transported from the final node of the supply chain to the customer's site has to report this malfunction in the supply chain operation. Assuming stochastic fluctuations driving the demand line, the production is adjusted based on the new value of demand; however the information mismatch is not taken into account. The integration of noise in the variables and not in the parameters models the problem in its fundamental base; it is the measurement of the information mismatch that needs to be modeled.

In this work, a two-stage supply chain is proposed where in the first stage, a MILP model is solved in order to provide the solutions for the second model; the levels of solutions that are derived from the first model concern the construction of the supply network and solutions that correspond to quantities transported throughout the supply chain. In the second model, the expected lead time is measured based on the amount of unsatisfied demand ($Δ$). Imposing thresholds on "small" or "large" $\Delta$, the network is reconstructed providing additional information regarding the capacity of the facilities and the magnitude of products that need to be constructed as it is assumed that warehouses serve as small production plants in order to minimize the expected lead time and therefore increase service level. Three types of statistical noise are examined; Normal (Gaussian), Lognormal and Pareto. The decision levels of the variables with Gaussian noise report larger fluctuation from the actual situation (deterministic), while the smaller fluctuation is observed from Pareto noise and especially with the smallest exponent. The cost function derived showcases how to incorporate such stochasticity in a supply chain model and what eventual benefits one may derive out of it. As a tailored example, we show that a producer may benefit from better return only through a suitable selective choice of producers whose production cost probability density function abides a Pareto distribution. Such a study can have a significant impact on any overall supply chain cost due to the linearly increasing objective function. While such stochastic optimization is not unknown in the realm of statistical mechanics Spall, (2003), the mapping is an altogether new concept in supply chain literature, an approach that has the prospect of coming up with rich dividends in the future.

**Appendix-Nomenclature**

**Index**

$i$  Plant

$j$  Warehouse

$k$  Customer

**Continuous Variables**

$TC$  Total supply chain cost

$P_i$  Produced quantities in plant $i$

$Q_{ij}$  Transported quantities from plant $i$ to warehouse $j$

$Q_{jk}$  Transported quantities from warehouse $j$ to customer $k$

$E_{jk}$  Transported quantities from warehouse $j$ to customer $k$ that exceed a certain level (high)

$\hat{P}_i$  Produced quantities in plant $i$ with noise representation

$\hat{Q}_{ij}$  Transported quantities from plant $i$ to warehouse $j$ with noise representation

$\hat{Q}_{jk}$ Transported quantities from warehouse $j$ to customer $k$ with noise representation

$R_{jk}$ Transported quantities from warehouse $j$ to customer $k$ that exceed a certain level (low)

$W_j$ Capacity of warehouse $j$

$ELD_k$ Expected Lead Time of customer $k$

$P_k(Q^U)$ Stock out probability of customer $k$

$P_k(Q^L)$ Overstocking probability of customer $k$

$\Delta_k$ Deficit in demand satisfaction for customer $k$

**Binary Variables**

$X_{ij}$ 1 if the corresponding connection between plant $i$ and warehouse $j$ exists, 0 otherwise

$X_{jk}$ 1 if the corresponding connection between warehouse $j$ and customer $k$ exists, 0 otherwise

$Y_j$ 1 if warehouse $j$ is selected, 0 otherwise

$K_{jk}$ 1 if small quantities will be delivered from warehouse $j$ to customer $k$ due to large demand deficit, 0 otherwise

$\Omega_{jk}$ 1 if large quantities will be delivered from warehouse $j$ to customer $k$ due to large demand deficit, 0 otherwise

$\lambda_k$ 1 if the deficit in demand satisfaction is lies in the interval $[\Delta^L, \Delta^U]$, 0 otherwise

$\zeta_k$ 1 if the deficit in demand satisfaction is lies in the interval $[\Delta^U, \Delta]$, 0 otherwise

**Parameters**

$P_i^U$ Upper bounded production of plant $i$

$P_i^L$ Lower bounded production of plant $i$

| | |
|---|---|
| $Q_{ij}^U$ | Maximum capacity of transported quantities from plant $i$ to warehouse $j$ |
| $Q_{ij}^L$ | Minimum capacity of transported from plant $i$ to warehouse $j$ |
| $Q_{jk}^U$ | Maximum capacity of transported from warehouse $j$ to customer $k$ |
| $Q_{jk}^L$ | Minimum capacity of transported from warehouse $j$ to customer $k$ |
| $I_j$ | Inventory held at warehouse $j$ |
| $D_k$ | Demand of customer $k$ |
| $\Delta_k$ | Stockout quantity in customer $k$ |
| $T^u$ | Maximum time for product delivery |
| $T^l$ | Minimum time for product delivery |
| $a_j$ | Coefficient relating quantity at capacity at warehouse j |
| $\beta_{jk}$ | Production rate for quantities stored at warehouse $j$ that will be delivered to customer $k$ in order to cover the high deficit in demand satisfaction. |
| $\gamma_{jk}$ | Production rate for quantities stored at warehouse $j$ that will be delivered to customer $k$ in order to cover the low deficit in demand satisfaction. |

**Cost parameters**

| | |
|---|---|
| $c_i^P$ | Production cost of plant $i$ |
| $c_{ij}^{VTR}$ | Variable transportation cost of plant $i$ to warehouse $j$ |
| $c_{ij}^{FTR}$ | Fixed transportation cost of plant $i$ to warehouse $j$ |
| $c_{jk}^{VTR}$ | Variable transportation cost of warehouse $j$ to customer $k$ |
| $c_{jk}^{FTR}$ | Fixed transportation cost of warehouse $j$ to customer $k$ |
| $c_j^{IN}$ | Installation cost of warehouse $j$ |
| $c_{jk}^{PO}$ | Production cost of small quantities that will be manufactured in warehouse $j$ and will be delivered to customer $k$ |
| $c_{jk}^{PU}$ | Production cost of large quantities that will be manufactured in warehouse $j$ and will be delivered to customer $k$ |

**References**

Acar, Yavuz, Sukran Kadipasaoglu, and Peter Schipperijn. 2009. "A Decision Support Framework for Global Supply Chain Modelling: An Assessment of the Impact of Demand, Supply and Lead-Time Uncertainties on Performance." *International Journal of Production Research* 48 (11): 3245–68. doi:10.1080/00207540902791769.

An, N., J.-C. Lu, D. Rosen, and L. Ruan. 2007. "Supply-Chain Oriented Robust Parameter Design." *International Journal of Production Research* 45 (23): 5465–84. doi:10.1080/00207540701325124.

Arabatzis, Garyfallos, Konstantinos Petridis, Spyros Galatsidas, and Konstantinos Ioannou. 2013. "A Demand Scenario Based Fuelwood Supply Chain: A Conceptual Model." *Renewable and Sustainable Energy Reviews* 25: 687–97.

Arabatzis, Garyfallos, Konstantinos Petridis, and P. Kougioulis. 2014. "Proposing a Supply Chain Model for the Production-Distribution of Fuelwood in Greece Using Multi-Objective Programming." *E-Innovation for Sustainable Development of Rural Resources during Global Economic Crisis, IGI Global,(forthcoming)*, 171–80.

Baghalian, Atefeh, Shabnam Rezapour, and Reza Zanjirani Farahani. 2013. "Robust Supply Chain Network Design with Service Level against Disruptions and Demand Uncertainties: A Real-Life Case." *European Journal of Operational Research* 227 (1): 199–215. doi:10.1016/j.ejor.2012.12.017.

Beamon, Benita M. 1998. "Supply Chain Design and Analysis::: Models and Methods." *International Journal of Production Economics* 55 (3): 281–94.

Fisher, Marshall L. 1985. "An Applications Oriented Guide to Lagrangian Relaxation." *Interfaces* 15 (2): 10–21.

Gebreslassie, Berhane H., Yuan Yao, and Fengqi You. 2012. "Design under Uncertainty of Hydrocarbon Biorefinery Supply Chains: Multiobjective Stochastic Programming Models, Decomposition Algorithm, and a Comparison between CVaR and Downside Risk." *AIChE Journal* 58 (7): 2155–79.

Grigoroudis, Evangelos, Konstantinos Petridis, and Garyfallos Arabatzis. 2014. "RDEA: A Recursive DEA Based Algorithm for the Optimal Design of Biomass Supply Chain Networks." *Renewable Energy* 71: 113–22.

Gruen, Thomas, and Daniel Corsten. 2004. "Stock-Outs Cause Walkouts." *Harvard Business Review* 82 (5): 26–27.

Gruen, T. W., D. S. Corsten, and S. Bharadwaj. 2002. "Retail Out-of-Stocks." *A Worldwide Examination of Extent, Causes, and Consumer Responses, Research Study, Atlanta*.

Guillén-Gosálbez, Gonzalo, and Ignacio E. Grossmann. 2009. "Optimal Design and Planning of Sustainable Chemical Supply Chains under Uncertainty." *AIChE Journal* 55 (1): 99–121.

Jindal, Anil, and Kuldip Singh Sangwan. 2014. "Closed Loop Supply Chain Network Design and Optimisation Using Fuzzy Mixed Integer Linear Programming Model." *International Journal of Production Research* 52 (14): 4156–73. doi:10.1080/00207543.2013.861948.

Kim, Jinkyung, Matthew J. Realff, and Jay H. Lee. 2011. "Optimal Design and Global Sensitivity Analysis of Biomass Supply Chain Networks for Biofuels under Uncertainty." *Computers & Chemical Engineering* 35 (9): 1738–51.

Krikke, Harold, Jacqueline Bloemhof-Ruwaard, and Luk N. Van Wassenhove. 2001. *Design of Closed Loop Supply Chains: A Production and Return Network for Refrigerators*. Erasmus Research Institute of Management (ERIM). https://flora.insead.edu/fichiersti_wp/inseadwp2001/2001-67.pdf.

Krikke, H., J. Bloemhof-Ruwaard, and L. N. Van Wassenhove. 2003. "Concurrent Product and Closed-Loop Supply Chain Design with an Application to Refrigerators." *International Journal of Production Research* 41 (16): 3689–3719.

Lee, Hau L., Venkata Padmanabhan, and Seungjin Whang. 2004. "Information Distortion in a Supply Chain: The Bullwhip Effect." *Management Science* 50 (12_supplement): 1875–86.

Liu, Songsong, and Lazaros G. Papageorgiou. 2013. "Multiobjective Optimisation of Production, Distribution and Capacity Planning of Global Supply Chains in the Process Industry." *Omega* 41 (2): 369–82.

Melo, M. Teresa, Stefan Nickel, and Francisco Saldanha da Gama. 2006. "Dynamic Multi-Commodity Capacitated Facility Location: A Mathematical Modeling Framework for Strategic Supply Chain Planning." *Computers & Operations Research* 33 (1): 181–208.

Mohammadi Bidhandi, Hadi, Rosnah Mohd. Yusuff, Megat Mohamad Hamdan Megat Ahmad, and Mohd Rizam Abu Bakar. 2009. "Development of a New Approach for Deterministic Supply Chain Network Design." *European Journal of Operational Research* 198 (1): 121–28. doi:10.1016/j.ejor.2008.07.034.

Pan, Feng, and Rakesh Nagi. 2010. "Robust Supply Chain Design under Uncertain Demand in Agile Manufacturing." *Computers & Operations Research* 37 (4): 668–83.

Petridis, Konstantinos. 2013. "Optimal Design of Multi-Echelon Supply Chain Networks under Normally Distributed Demand." *Annals of Operations Research*, July, 1–29. doi:10.1007/s10479-013-1420-6.

Riddalls, C. E., and S. Bennett. 2002. "Production-Inventory System Controller Design and Supply Chain Dynamics." *International Journal of Systems Science* 33 (3): 181–95. doi:10.1080/00207720110092180.

Santoso, Tjendera, Shabbir Ahmed, Marc Goetschalckx, and Alexander Shapiro. 2005. "A Stochastic Programming Approach for Supply Chain Network Design under Uncertainty." *European Journal of Operational Research* 167 (1): 96–115.

Seferlis, Panos, and Nikolaos F. Giannelos. 2004. "A Two-Layered Optimisation-Based Control Strategy for Multi-Echelon Supply Chain Networks." *Computers & Chemical Engineering* 28 (5): 799–809.

Spearman, Mark L., and Michael A. Zazanis. 1992. "Push and Pull Production Systems: Issues and Comparisons." *Operations Research* 40 (3): 521–32.

Tamas, Mick. 2000. "Mismatched Strategies: The Weak Link in the Supply Chain?" *Supply Chain Management: An International Journal* 5 (4): 171–75.

Tsao, Yu-Chung, and Jye-Chyi Lu. 2012. "A Supply Chain Network Design Considering Transportation Cost Discounts." *Transportation Research Part E: Logistics and Transportation Review* 48 (2): 401–14.

Tsiakis, Panagiotis, Nilay Shah, and Constantinos C. Pantelides. 2001. "Design of Multi-Echelon Supply Chain Networks under Demand Uncertainty." *Industrial & Engineering Chemistry Research* 40 (16): 3585–3604.

Wang, Fan, Xiaofan Lai, and Ning Shi. 2011. "A Multi-Objective Optimization for Green Supply Chain Network Design." *Decision Support Systems* 51 (2): 262–69.

You, Fengqi, and Ignacio E. Grossmann. 2008. "Design of Responsive Supply Chains under Demand Uncertainty." *Computers & Chemical Engineering* 32 (12): 3090–3111.